\documentclass[9pt,twocolumn,twoside]{pnas-new}
\usepackage[linesnumbered,ruled,vlined,algo2e]{algorithm2e}
\usepackage{amsthm}
\usepackage{cleveref}
\usepackage{amsmath}
\usepackage{multirow}
\usepackage{graphicx}
\usepackage{subcaption}
\usepackage{pdfpages}

\newcommand{\bm}[1]{{\boldsymbol{#1}}}
\newcommand{\SI}[1]{{\textcolor{blue}{#1}}}

\newcommand{\appDetJ}[1]{{\det \mathcal{J}}_{#1}}

\newcommand{\assign}{\leftarrow}
\SetAlgorithmName{Pseudocode}{}{}

\newcounter{algoCount}
\newcommand{\pc}[2]{
\quad\quad 
\begin{minipage}{.8\linewidth}
\setcounter{algorithm}{1}
\begin{algorithm}[H]
\renewcommand{\thealgorithm}{\thealgoCount}
\caption{#1}
#2
\end{algorithm}
\end{minipage}
\addtocounter{algoCount}{1}
}

\DeclareCaptionFormat{overlay}{\gdef\capoverlay{#1#2#3\par}}
\DeclareCaptionStyle{overlay}{format=overlay}

\templatetype{pnasresearcharticle} 

\newtheorem{theorem}{Theorem}
\newtheorem{corollary}{Corollary}

\title{Improving sampling efficacy on high dimensional distributions with thin high density regions using Conservative Hamiltonian Monte Carlo}

\author[a,1]{Geoffrey McGregor}
\author[b,1,2]{Andy T.S. Wan} 

\affil[a]{Department of Mathematics, University of Toronto, ON M5S2E4}
\affil[b]{Department of Applied Mathematics, University of California, Merced, CA 95343}

\leadauthor{McGregor} 

\authorcontributions{G.M. and A.T.S.W. performed research, analysis, computation and wrote the paper together.}
\authordeclaration{The authors declare no conflict of interest.}
\equalauthors{\textsuperscript{1} G.M. and A.T.S.W. contributed equally to this work.}
\correspondingauthor{\textsuperscript{2}To whom correspondence should be addressed. E-mail: \href{mailto://andywan@ucmerced.edu}{andywan@ucmerced.edu}}

\keywords{Markov Chain Monte Carlo $|$ Hamiltonian Monte Carlo $|$ energy-preserving integrator $|$ approximate stationarity} 

\begin{abstract}
Hamiltonian Monte Carlo is a prominent Markov Chain Monte Carlo algorithm, which employs symplectic integrators to sample from high dimensional target distributions in many applications, such as statistical mechanics, Bayesian statistics and generative models. However, such distributions tend to have thin high density regions, posing a significant challenge for symplectic integrators to maintain the small energy errors needed for a high acceptance probability. Instead, we propose a variant called Conservative Hamiltonian Monte Carlo, using $R$--reversible energy-preserving integrators to retain a high acceptance probability. We show our algorithm can achieve approximate stationarity with an error determined by the Jacobian approximation of the energy-preserving proposal map. Numerical evidence shows improved convergence and robustness over integration parameters on target distributions with thin high density regions and in high dimensions. Moreover, a version of our algorithm can also be applied to target distributions without gradient information.
\end{abstract}

\dates{This manuscript was compiled on \today}

\begin{document}
\maketitle
\ifthenelse{\boolean{shortarticle}}{\ifthenelse{\boolean{singlecolumn}}{\abscontentformatted}{\abscontent}}{}
For more than half a century, Markov Chain Monte Carlo (MCMC) algorithms have been utilized in numerous applications across science and engineering, from its early days in statistical mechanics \cite{metropolis1953equation, hastings1970monte} to Bayesian statistics \cite{Tierney94,gelman1995bayesian}, and more recently in generative models \cite{SongAE21,YangAE23}. A gradient-based MCMC algorithm known as Hamiltonian Monte Carlo (HMC) \cite{duane1987hybrid, neal1994improved, neal2011mcmc, betancourt2017conceptual} has seen recent wide adoption for many applications in Bayesian statistics. Specifically, given a target distribution $\pi(\bm \theta)$, HMC extends the sample space by interpreting $\bm \theta$ as generalized coordinate variables $\bm q\in \mathbb{R}^d$ and introducing momentum variables $\bm p\in \mathbb{R}^d$ giving rise to a joint distribution $\pi(\bm q, \bm p) \propto \exp(-H(\bm q,\bm p))$, where $H(\bm q,\bm p) = K(\bm p) + U(\bm q)$ is the associated Hamiltonian function with the kinetic energy $K(\bm p) = \frac{1}{2}\bm p^T M^{-1} \bm p$ and the potential energy $U(\bm q) = -\log \pi(\bm q)$. From a current sample $(\bm q^i, \bm p^i)$, a new proposal $ (\bm q^*,\bm p^*)$ is obtained by numerically solving the associated Hamiltonian system $\dot{\bm q} = M^{-1}\bm p, \dot{\bm p} = -\nabla U(\bm q)$ over a prescribed time interval $t\in [0,T]$ using a symplectic integrator \cite{hair06Ay} of step size $\tau$, where typically a Leapfrog or Str\"omer-Verlet integrator is used. Under appropriate conditions, HMC would satisfy the stationarity condition \cite{neal2011mcmc}, ensuring the generated samples converging to the target distribution.
\\\indent One advantage of employing such a proposal map is far distant samples can be obtained via Hamiltonian dynamics, thus improving sampling efficiency over traditional random-walk MCMC algorithms. Moreover, as symplectic integrators preserve volume (i.e. $\det J_{\Psi_{SYM}}=1$ where $\Psi_{SYM}$ denotes a symplectic proposal map), the Metropolis adjustment step can avoid the costly computation of the Jacobian of the proposal map $J_{\Psi_{SYM}}$ \cite{neal2011mcmc}, leading to HMC having an acceptance probability of $\alpha_{SYM} = \min\left(1,\exp(-\Delta H)\right)$ where $\Delta H$ is the energy difference between the proposed and current samples. For the Leapfrog integrator, it is known that $\Delta H = \mathcal{O}(\tau^2)$ via backward error analysis \cite{hair06Ay,beskos2013optimal}, leading to a high acceptance probability if $\tau$ is sufficiently small. Moreover, progress has been made to tune the integration parameters $\tau, T$ and mass matrix $M$, such as No--U--Turn sampling \cite{hoffman2014no}, tuning step sizes \cite{beskos2013optimal} and generalizing to $M(\bm q)$ in Riemannian HMC \cite{girolami2011riemann, betancourt2017conceptual}.
\\\indent However, despite the successes of HMC, there remains aspects of the algorithm which can still be improved. In particular, symplectic integrators do not in general preserve the Hamiltonian exactly. As the acceptance probability $\alpha$ depends on the error in the Hamiltonian, this can lead to more rejected proposals as the dimension $d$ increases. Indeed, as discussed in \cite{beskos2013optimal} with suitable regularity assumptions on $\pi$, the step size $\tau$ used in the Leapfrog integrator of HMC must scale as $\mathcal{O}(d^{-\frac{1}{4}})$, in order to maintain a constant acceptance probability as $d$ increases. One intuitive explanation behind this performance decrease is that high dimensional distributions can concentrate on thin high density regions \cite{vershynin18}. 
     Thus, increasing $d$ leads to a decrease in sampling efficacy for HMC, as its proposals are likely to be rejected with Leapfrog integrator being unable to remain near the constant energy surface of thin high density regions. Instead, we propose using energy-preserving integrators to allievate this difficulty in sampling from high dimensional distributions with concentrated high density regions.
\section{Conservative Hamiltonian Monte Carlo (CHMC)}
In order to obtain samples which stay on the same Hamiltonian or energy level set (up to machine precision) after numerical integration, we propose to use energy-preserving integrators, instead of symplectic integrators\footnote{There are no known general integrators which can simultaneously preserve energy and be symplectic, as such integrator would be equivalent to a time-reparametization of exact solutions \cite{mars88Ay}.}. From the field of geometric numerical integration \cite{hair06Ay}, there are a number of well-known energy-preserving integrators\footnote{There are also other approaches which preserve energy, such as projection methods \cite{hair06Ay} and relaxation methods \cite{dekkerVerwer84,CalvoAE06} but they do not in general satisfy $R$--reversibility.}, such as the Itoh--Abe Discrete Gradient scheme \cite{IA88}, Average Vector Field (AVF) Discrete Gradient scheme \cite{quis08a}, or Discrete Multiplier Method (DMM) \cite{WBN2017}. Employing any of these approaches could be used within our proposed algorithm, called \emph{Conservative Hamiltonian Monte Carlo} (CHMC). 
\hskip -4mm\pc{CHMC Algorithm}{
	Pick $\bm \theta^0$
	\For{$i = 1,2,\dots, K$}{
		$\bm q^0 \assign \bm \theta^i $\\
		Draw $\bm p^0 \sim \mathcal{N}(\bm 0, M)$\\
		$(\bm q^*, \bm p^*) \assign \Psi_{EP}(\bm q^0, \bm p^0; H, \tau, T)$\\
        $\Delta H \assign H(\bm q^*, \bm p^*)-H(\bm q^0,\bm p^0)$\\
		$\alpha \assign \min\left\{1,\exp\left(-\Delta H\right)\appDetJ{\Psi_{EP}}(\bm q^0,\bm p^0)\right\}$\\
		$ \bm \theta^{i+1} \assign \begin{cases}&\bm q^*\text{ with probability } \alpha\\ & \bm \theta^i, \text{ with probability } 1-\alpha \end{cases}$
	}\label{CHMCAlg}
}\\

The CHMC algorithm is similar to that of HMC but with two distinct differences. First is the usage of an energy-preserving proposal map $\Psi_{EP}$, where a new proposal is obtained by integrating the associated Hamiltonian system using an energy-preserving integrator, starting at $(\bm q^0,\bm p^0)$ with 
a uniform step size $\tau$ and a fixed integration length $T$ \footnote{Adaptive step size and integration length can also be used, such as with No--U--Turn sampling \cite{hoffman2014no}.}. In contrast, HMC typically utilizes symplectic integrators, or other integrators which do not preserve energy \cite{fang2014compressible}. Second is the appearance of the Jacobian approximation\footnote{ $\displaystyle\appDetJ{\Psi_{EP}}$ can include dependence on $(\bm q^*, \bm p^*)$, defined implicitly through $\Psi_{EP}(\bm q^0,\bm p^0)$.} of the proposal map, $\appDetJ{\Psi_{EP}}$, due to the non-volume preserving transformation of such energy-preserving integrators. While exact stationarity can be achieved for CHMC using the full Jacobian, as discussed in \SI{\emph{SI Appendix}, Section A}, we instead propose using an approximate Jacobian to strike a balance between computational efficiency and approximate stationarity. Our next theorem states \Cref{CHMCAlg} with an approximate Jacobian $\appDetJ{\Psi}$ achieving \emph{approximate stationarity} to the target distribution, as proved in \SI{\emph{SI Appendix}, Section A}.
\begin{theorem}[Error bound on stationarity of $R$--reversible proposal with approximate Jacobian]\label{CHMCStat} Denote $\bm z :=(\bm{q},\bm{p}) \in\mathbb{R}^{2d}$ and let $\Psi:\mathbb{R}^{2d}\rightarrow \mathbb{R}^{2d}$ be a positively-oriented (i.e. $\det J_\Psi > 0$) $C^1$-diffeomorphism, with its Jacobian matrix entries $[J_\Psi]_{ij} \in  L^{\infty}(\mathbb{R}^{2d})$. Also, suppose $\Psi$ is $R$--reversible \cite{hair06Ay} with respect to the bijection $R(\bm{q},\bm{p})=(\bm{q},-\bm{p})$, i.e. $R\circ\Psi\circ R\circ\Psi=I$, and let $\appDetJ{\Psi}\in L^{\infty}(\mathbb{R}^{2d})$ be an approximation of $\det J_\Psi$. Denoting the error $\epsilon(\bm z):=\det J_{\Psi}(\bm z)-\appDetJ{\Psi}(\bm z) \in L^{\infty}(\mathbb{R}^{2d})$ and letting $\pi \in L^1(\mathbb{R}^{2d})\cap L^{\infty}(\mathbb{R}^{2d})$ be a target density satisfying $\pi\circ R=\pi$, define the acceptance probability to be 
\begin{equation*}
    \alpha({\bm z}) = \min\left(1 ,\frac{\pi\left(\Psi({\bm z})\right)}{\pi({\bm z})}\appDetJ{\Psi}({\bm z})\right), 
\end{equation*}
with the transition kernel density from $\bm z$ to $\bm z'$ be given by $\rho({\bm z},{\bm z}')=\alpha({\bm z})\delta({\bm z}'-\Psi({\bm z}))+(1-\alpha({\bm z}))\delta({\bm z}'-R({\bm z}))$ \cite{fang2014compressible, bou2018geometric}, where $\delta(\cdot)$ denotes the Dirac distribution in $\mathbb{R}^{2d}$. Then, the error from stationarity can be bounded as
\begin{equation*}
\biggr|\pi({\bm z}')-\int\limits_{\mathbb{R}^{2d}}\rho({\bm z},{\bm z}')\pi({\bm z})\text{d}{\bm z}\biggr| \leq C(\bm z') \pi(\bm z') \|\epsilon\|_\infty,
\end{equation*}where $C(\bm z')=\frac{}{}2+\mathcal{O}(\Delta H)+\left|\frac{}{}1-\det J_{\Psi}({R(\bm z')})\right|$.
\end{theorem}
\noindent The next corollary shows a stationarity-error bound for proposal maps of energy-preserving integrators and a lower bound on acceptance probability, see \SI{\emph{SI Appendix}, Section B and C}.
\begin{corollary}[Approximate Stationarity of CHMC]\label{cor: CHMCStat}
Let $\Psi_{EP}$ satisfy the hypotheses of \Cref{CHMCStat} and be an $N$-times composition of an energy-preserving integrator with a uniform step size $\tau$ such that $\det J_{\Psi_{EP}}(\bm z)=1+C(\bm z)\tau^p+O(\tau^{2p})$, where $C\in L^{\infty}(\mathbb{R}^{2d})$ for some $p>0$. Then for $\appDetJ{\Psi_{EP}}(\bm z)=1$ and sufficiently small $\tau$, \Cref{CHMCAlg} satisfies the approximate stationarity result of \Cref{CHMCStat} with an error of $\|\epsilon\|_\infty = \mathcal{O}(\tau^p)$ and the acceptance probability is bounded below by $e^{-\delta}$, for any desired energy error tolerance $\delta>0$. 
\end{corollary}
\subsection*{Implementation details}
For a general target distribution $\pi$, numerical schemes which preserve $H$ are typically implicit, where a nonlinear system needs to be solved at each step $\tau$ using iterative methods such as fixed point iterations, quasi-Newton or Newton methods. For sufficiently small step size, each iteration reduces the residual of the current energy error $\Delta H$ until it reaches below a desired energy tolerance $\delta$. Thus, for CHMC to be efficient in practice, a balance needs to be struck between the energy tolerance $\delta$ and the number of iterations to solve the implicit energy-preserving scheme.
\\\indent As discussed in \Cref{cor: CHMCStat} for $p=2$, employing a second-order energy-preserving $R$--reversible scheme in \Cref{CHMCAlg} with $\appDetJ{\Psi_{EP}}=1$ leads to samples satisfying approximate stationarity with an error of $\mathcal{O}(\tau^2)$. For instance, the symmetrized Itoh--Abe Discrete Gradient or DMM scheme (\SI{\emph{SI Appendix}, Section E.1}), or the AVF scheme (\SI{\emph{SI Appendix}, Section E.2}), are second-order energy-preserving $R$--reversible schemes, as detailed in \SI{\emph{SI Appendix}, Sections F--H}. Specifically, the AVF scheme requires gradient information of $H$, followed typically by a quadrature approximation of an integral associated with the scheme. In contrast, the symmetrized scheme of the Itoh--Abe Discrete Gradient or DMM does not require gradient information of $H$, but with potential regularization needed for small divisors. 
\\\indent Moreover, the error from stationarity can be further reduced by choosing $\appDetJ{\Psi_{EP}}$ to be a higher order approximation of $\det J_{\Psi_{EP}}$, such as using an improved approximation of the determinant involving traces, as detailed in \SI{\emph{SI Appendix}, Section I}. 
However, higher order approximations of the determinant generally add computational costs, which may outweigh the benefits of the improved error from stationarity. In our numerical experiments, we have employed CHMC with $\appDetJ{\Psi_{EP}}=1$, to minimize the computational cost of the Jacobian. It's also worthwhile to point out that choosing $\appDetJ{\Psi_{EP}}=1$ and an energy-preserving $R$--reversible integrator not requiring gradient information of $H$ yields a $H$--\emph{gradient--free} version of \Cref{CHMCAlg}, such as using the symmetrized Itoh--Abe Discrete Gradient scheme or symmetrized DMM scheme. The gradient--free CHMC may be useful in applications where first derivative of the target distribution is not readily accessible.
\section{Results}
We compare the sampling efficacy of CHMC versus HMC on two target distributions, $p$-generalized $\chi$ distribution and $p$-generalized Gaussian distribution \cite{richter07}. Unless stated otherwise, we compare HMC with the Leapfrog integrator (HMC--LF) and CHMC with the symmetrized Itoh--Abe Discrete Gradient or DMM scheme, with a maximum of two fixed-point iterations in Newton's method, and an approximate Jacobian $\appDetJ{\Psi_{EP}}=1$, as detailed in \SI{\emph{SI Appendix}, Section J}. Also, for a baseline comparison, a uniform step size $\tau$ and a fixed integration length $T$ are used for both algorithms, with adaptive step size and variable integration length left for future work.
\subsection*{$p$-generalized $\chi$ distribution}
\captionsetup[figure]{font=small,skip=0pt}
\begin{figure*}[ht!]
    \centering        
    \includegraphics[scale=.48,trim = 2.1cm 0.6cm 0cm 0.4cm, clip]{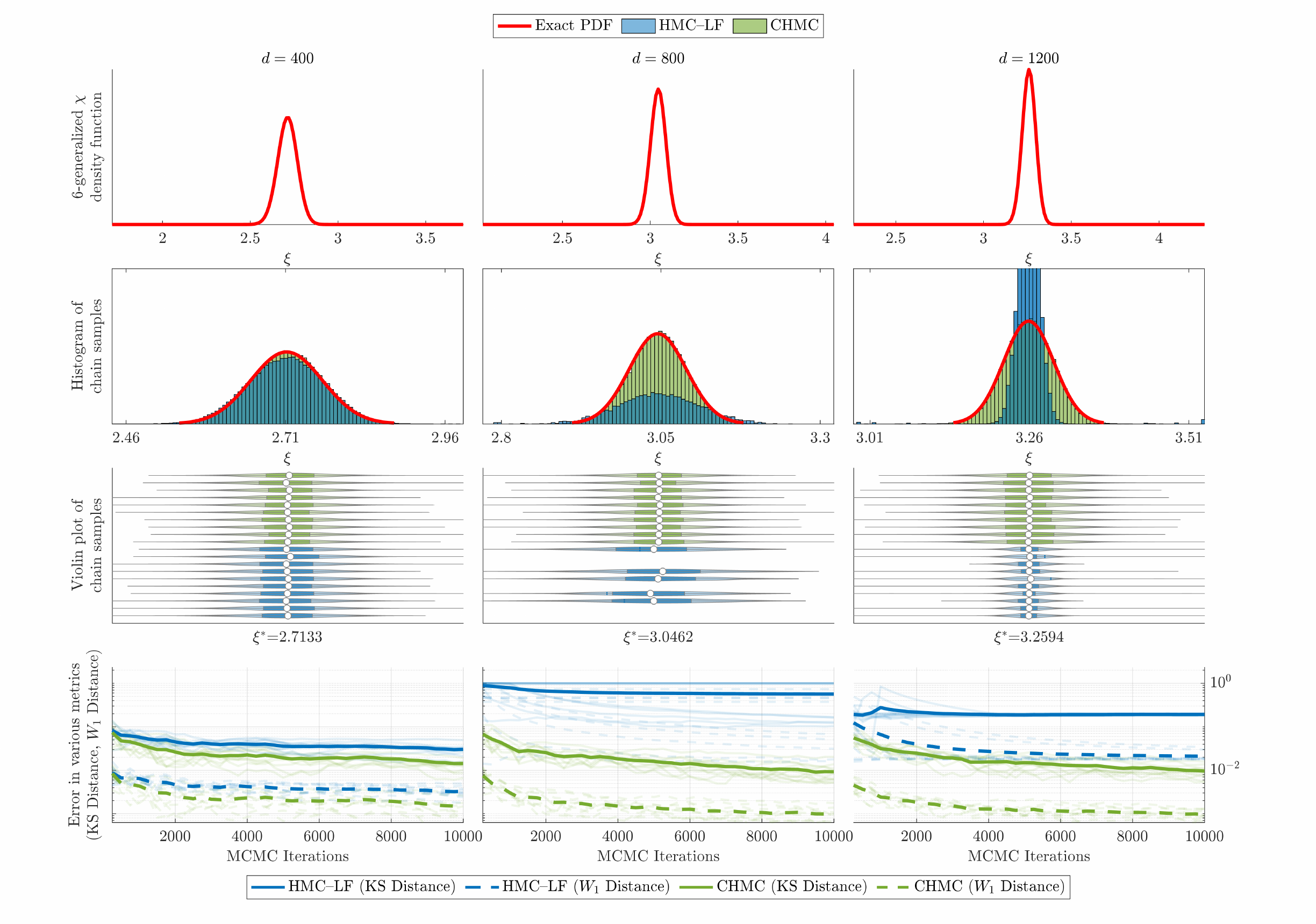}
    \caption{Comparison on histograms and convergence of HMC--LF versus CHMC at sampling the $6$-generalized $\chi$ distribution with increasing degrees of freedom $d$.
    }\label{fig:Hist1} 
\end{figure*}
\begin{figure*}[ht!]
    \centering
    \begin{subfigure}[t]{\textwidth}
        \centering
        \includegraphics[scale=.53,trim = 1.5cm 0.5cm 0cm 0cm, clip]{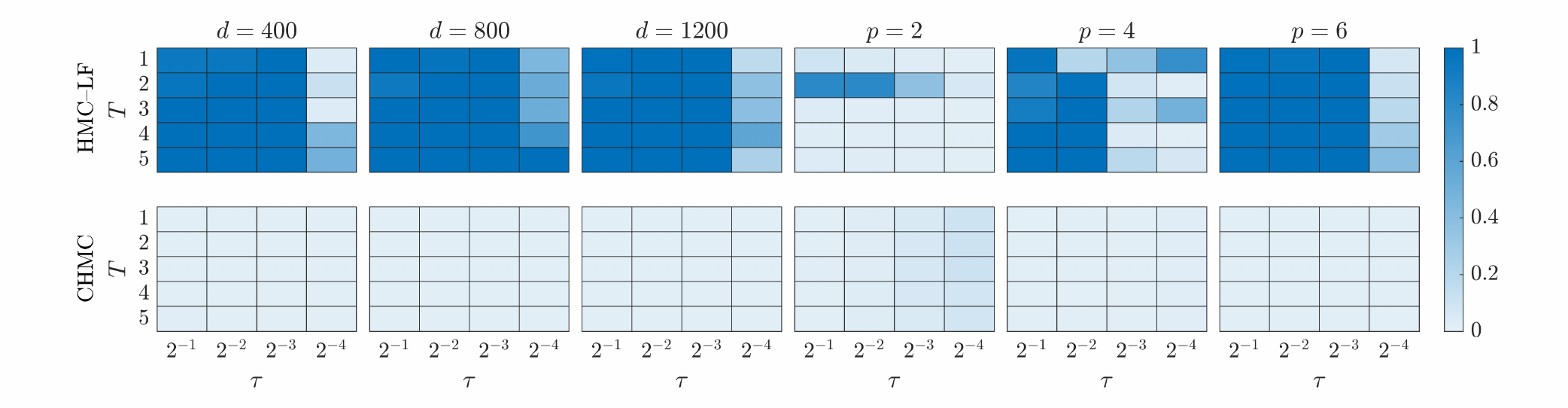}
    \caption{Heat map of errors in Kolmogorov--Smirnoff distance.}
    \label{fig:KSplot}
    \end{subfigure}
    \centering
    \begin{subfigure}[t]{\textwidth}
        \centering
        \includegraphics[scale=.53,trim = 1.5cm 0.5cm 0cm 0cm, clip]{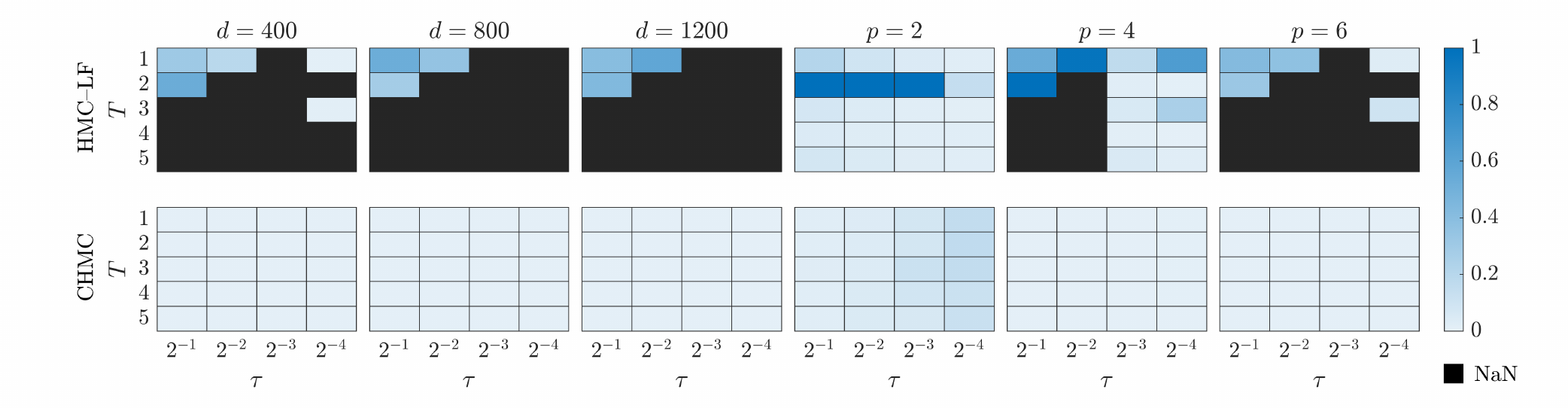}
    \caption{Heat map of errors in Wasserstein $W_1$ distance.}
    \label{fig:W1plot}
    \end{subfigure}
    \caption{Comparison on errors of HMC--LF versus CHMC at sampling $p$-generalized $\chi$ distribution for various $d, p$, integration parameters $T,\tau$.}
\end{figure*}
We first demonstrate the sampling efficacy of CHMC over HMC on distributions with concentrated high density regions. Specifically, the $p$-generalized $\chi$-distribution with $d$ degrees of freedom \cite{richter07} has the density function
$\pi_{d,p}(\xi) = I_{(0,\infty)}(\xi)p^{1-\frac{d}{p}}\Gamma\left(\frac{d}{p}\right)^{-1}\xi^{d-1}\exp\left(-\frac{\xi^p}{p}\right).$
\\\indent As shown in Lemma 1 of \SI{\emph{SI Appendix}, Section D}, the majority of the density of $\pi_{d,p}$ lies in an annulus in $\ell_p$ norm with its width decreasing exponentially when $p>2$ and $d\gg 1$. 
Specifically, increasing $p$ while fixing $d$, or vice-versa, increasing $d$ while fixing $p>2$, leads to an exponential decrease in the width of the interval $\mathcal{I}_{d,p}(\epsilon)$ centred at $\xi^*=(d-1)^{\frac{1}{p}}$ which contains the majority of the density. Thus, despite being a 1-dimensional distribution, the exponentially thinning of the density $\pi_{d,p}$ makes this an ideal test case for comparison.
\\\indent The first set of results are shown in \Cref{fig:Hist1}, where the sampling efficacy of HMC--LF and CHMC are compared as the degrees of freedom $d$ is varied. The first row of panels illustrate the thinning of the probability density $\pi_{d,p}(\xi)$ as $d$ is increased. The second row shows the histograms of combined samples generated by HMC--LF in blue, and CHMC in green, with the third row showing the associated violin plots across all ten chains. These results highlight HMC--LF's slower rate of convergence due to the thin density region, as well as non-convergence due to instability of the Leapfrog integrator. In contrast, CHMC continues to sample the target distribution effectively using the energy-preserving integrator. 
Additionally as detailed in \SI{\emph{SI Appendix}, Section K.3}, the bottom row of \Cref{fig:Hist1}  highlights the improved convergence of CHMC versus HMC--LF measured in the Wasserstein--1 \cite{vallender74} and Kolmogorov--Smirnov \cite{kolmogorov33,smirnov39} distances, as the number of MCMC iteration increases. 
For details on iterations, chains, and integration parameters, see \SI{\emph{SI Appendix}, Section K}. 
\\\indent Next we illustrate the robustness of these results across various parameters. \Cref{fig:KSplot} and \ref{fig:W1plot} show two sets of heat maps, comparing the sampling efficacy of HMC--LF with CHMC on two different metrics, by varying integration length $T$ and step size $\tau$, and the parameters $d$ and $p$. We first compare the errors in the Kolmogorov--Smirnoff distance, shown in \Cref{fig:KSplot}, and then using the Wasserstein--1 distance, shown in \Cref{fig:W1plot}. As seen in each heat map, HMC--LF is only able to sample the target distribution effectively when the step size $\tau$ is sufficiently small, leading to a decrease in sampling efficacy compared to CHMC. These results also show the sampling efficacy of HMC--LF is more sensitive to the integration parameters than CHMC, with CHMC yielding more consistent results over a wide range of integration parameters $T$ and $\tau$, and $d$ and $p$ values, as discussed in \SI{\emph{SI Appendix}, Section K.2}.
\subsection*{High dimensional I.I.D. $p$-generalized Gaussian}
We consider the family of independent identically distributed (I.I.D.) $p$-generalized Gaussian in $\mathbb{R}^d$ with the joint density
$\pi({\bm x}) \sim \exp \left(-p^{-1}\lVert{\bm x}\rVert_p^p\right)$.
Here, $1<p<\infty$ and $\lVert{\bm x}\rVert_p$ denotes the $\ell_p$ norm of a random vector ${\bm x} \sim \mathbb{R}^d$. Recalling from \cite[Theorem 6]{richter07}, the random variable $\xi:= \lVert{\bm x}\rVert_p$
is equivalent to the $p$-generalized $\chi$ distribution.   The main result of Lemma 1 in \SI{\emph{SI Appendix}, Section D} shows that the $p$-generalized Gaussian distribution has the majority of its density living on a thin-strip in $d$-dimensional space. In particular, as we increase the dimension $d$, for $p>2$, the width of this strip exponentially decreases, and therefore we expect HMC--LF's performance to decrease due to thinning of the high-density region.  
\begin{figure*}[ht!]
    \centering
    \hskip -2mm
    \includegraphics[scale=.47,trim = 0.5cm 0cm 0cm 0cm, clip]{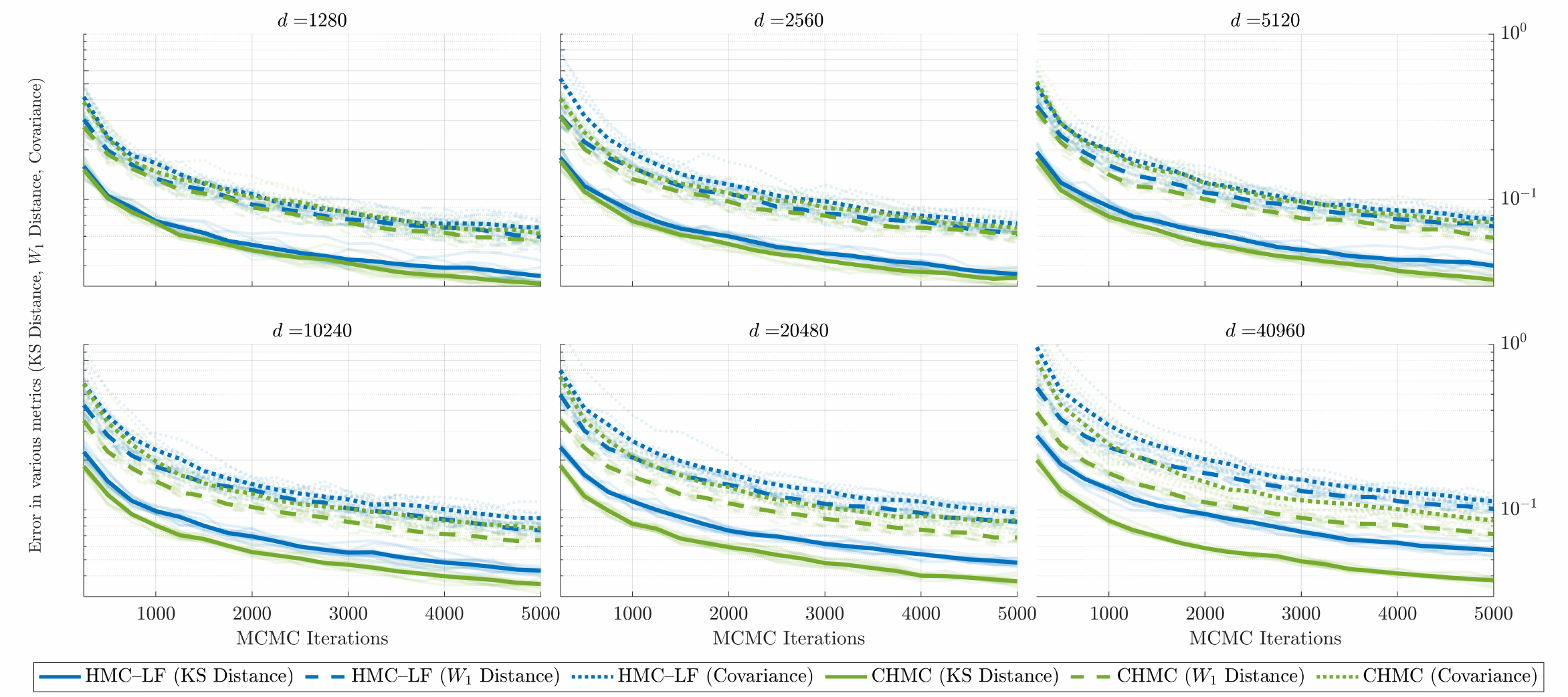}
    \caption{Comparison on convergence of HMC--LF versus CHMC in various metrics at sampling I.I.D. 4-generalized Gaussian in high dimensions.}
    \label{fig:convPGauss}
\end{figure*}
\\\indent \Cref{fig:convPGauss} consists of six convergence plots of increasing dimension $d$, each showing reduction in the errors as the number of MCMC iteration increases, measured in the Kolmogorov--Smirnov and Wasserstein--1 distances, as well as covariance. Since each component of $\bm x$ from the joint $p$-generalized Gaussian distribution is I.I.D., we computed the maximum of the two distances across each of their individual $d$ marginal distributions to save computational costs. Moreover, the error in covariance is also simplified and computed by taking the $l^{\infty}$ norm along the diagonal of the sample covariance matrix. As we observed with the $p-$generalized $\chi$ distribution, \Cref{fig:convPGauss} shows distinct separations measured in these metrics between HMC--LF and CHMC, as the width of the high-density region decreases. These results further highlight sampling performance in high-dimensional distributions with thin high-density regions can be improved by employing energy-preserving integrators. See \SI{\emph{SI Appendix}, Section L.1} for  details on \Cref{fig:convPGauss}.
\subsection*{Approximate stationarity and dimensional scaling of acceptance probability}
\begin{figure*}[ht!]
    \centering
    \includegraphics[scale=.51,trim = 2.5cm 0.5cm 0cm 1cm, clip]{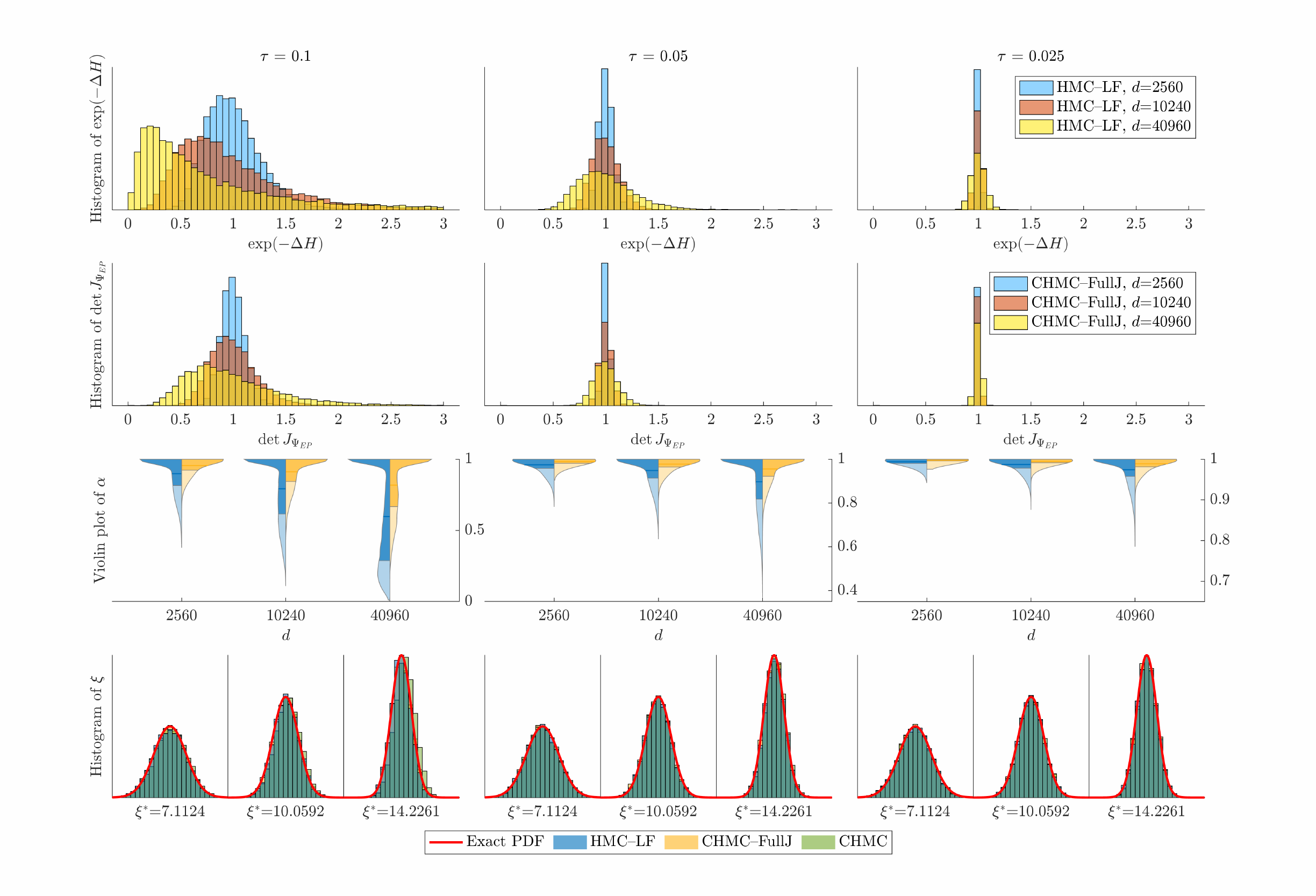}
    \caption{Comparison on histograms of $\exp(-\Delta H)$ versus $\det J_{\Psi_{EP}}$, violin plots of $\alpha$, and histogram of transformed samples (HMC--LF, CHMC--FullJ, CHMC).}
    \label{fig:compactComp}
\end{figure*}
So far, we have focused on CHMC with $\det \mathcal{J}_{\Psi_{EP}}=1$, with an acceptance probability of $\alpha_{EP} \geq \exp(-\delta)$. To observe the effects of approximate stationarity of CHMC in high dimensions, we first compare the acceptance probability of HMC--LF versus CHMC--FullJ (i.e. $\appDetJ{\Psi_{EP}}=\det J_{\Psi_{EP}}$), which satisfies exact stationarity as shown in \SI{\emph{SI Appendix}, Section A}. Specifically, the improvement on acceptance probability of CHMC--FullJ's $\alpha_{EP} = \min(1,\exp(-\delta) \det J_{\Psi_{EP}})$ over HMC--LF's $\alpha_{SYM} = \min(1,\exp(-\Delta H))$ hinges on favorable dimensional scaling of the Jacobian for conservative integrators over the negative exponential of the energy error for symplectic integrators.
\\\indent
The first and second rows in Figure \ref{fig:compactComp} each include three sets of histograms, comparing $\exp(-\Delta H)$ of HMC--LF versus $\det J_{\Psi_{EP}}$ of CHMC--FullJ, with increasing dimensions ($d=2560, 10240, 40960$) and decreasing step sizes ($\tau = 0.1, 0.05, 0.025$) across the columns. The third row of subfigures in \Cref{fig:compactComp} shows split violin plots comparing the acceptance probability of HMC--LF and CHMC--FullJ across the same dimensions and step sizes as above. As the first two rows of histograms illustrate, HMF-LF has larger variances on $\exp(-\Delta H)$ than CHMC--FullJ's $\det J_{\Psi_{EP}}$, across all $d$ and $\tau$. As a result, the violin plots show CHMC--FullJ has higher acceptance probability over HMC--LF, especially for large step sizes $\tau$. This indicates the growth of the Jacobian of the symmetrized Itoh--Abe scheme or symmetrized DMM scheme is slower than the growth of the energy error of the Leapfrog scheme, as $d$ increases. Moreover, since $\delta\approx 0$ for energy-preserving integrators, $\alpha_{EP}$ will concentrate near $1$ for sufficiently small $\tau$, which supports the choice of CHMC with $\appDetJ{\Psi_{EP}}=1$ as a reasonable approximation in practice.
\\\indent Finally, we look at the impacts of approximate stationarity in high dimensions by assessing at the effects of the transformation $\xi=\lVert\bm x\rVert_p$ on samples between the two equivalent distributions, the $p$-generalized Gaussian and $p$-generalized $\chi$ distributions \cite{richter07}. Specifically, we compare the transformed samples $\xi$ obtained from directly sampling the $p$-generalized Gaussian by HMC--LF, CHMC, and CHMC--FullJ. Despite CHMC's improvements over HMC--LF across various metrics discussed in previous examples, we do observe CHMC having a slight bias toward larger $\xi$ values as $d$ increases, corresponding to step size $\tau=0.1$ in the bottom left subfigure of \Cref{fig:compactComp}. One explanation for this bias is due to approximate stationarity of CHMC with the approximate Jacobian $\appDetJ{\Psi_{EP}}=1$, as outlined in Theorem \ref{CHMCStat}. Specifically, as indicated in the second row of subfigures in Figure \ref{fig:compactComp}, the exact Jacobian $\det J_{\Psi_{EP}} $ becomes less concentrated at 1 as $d$ increases, leading to extraneous samples being accepted by CHMC. In contrast, this bias is not present for CHMC--FullJ as it satisfies exact stationarity. On the other hand, since $\det J_{\Psi_{EP}}=1+\mathcal{O}(\tau^2)$  as shown in \SI{\emph{SI Appendix}, Section L.2}, reducing $\tau$ leads to $\det J_{\Psi_{EP}}$ being more concentrated at $1$, as seen in the second row of subfigures in Figure \ref{fig:compactComp}.  Thus by reducing $\tau$, the samples of CHMC and CHMC--FullJ become more similar, mitigating the observed bias without the need to compute the full Jacobian.

\section{Discussion}
We have introduced a variant of HMC, called CHMC, where an $R$--reversible energy-preserving integrator can be used to increase the acceptance probability and improve sampling efficacy of distributions with thin high density regions. To avoid computing the full Jacobian expression in the acceptance probability, an approximate Jacobian was introduced leading to the notion of approximate stationarity, where the associated error is determined by the choice of energy-preserving integrator, Jacobian approximation and step size. Our numerical studies showed various improvements of CHMC over HMC on the $p$-generalized $\chi$ and $p$-generalized Gaussian distributions across various parameters values and in high dimensions.

With the promising results of CHMC presented so far, there are various directions which warrant further investigation. For instance, alternative energy-preserving integrators can be explored for improved robustness and efficiency. Specifically, due to the implicit nature of $R$--reversible energy-preserving integrators discussed so far, more efficient nonlinear solvers can be investigated to improve computational costs. In addition, CHMC with adaptive step size and variable integration length can be explored, such as using No--U--Turn sampling. Also, it is of practical interest to apply CHMC and assess its effectiveness to large-scale applications, such as in statistical physics, Bayesian statistics and generative models. Moreover, the gradient--free aspect of the symmetrized Itoh--Abe or DMM scheme provides a promising alternative for HMC in applications with target distributions lacking derivative information. Additionally, a convergence theory for CHMC can be developed to assess how approximate stationarity influences potential bias and the convergence rate of CHMC.

\section{Materials and methods}

The implementation details are described in \SI{\emph{SI Appendix}, Sections K--L}. The MATLAB codes are available at the repository: \href{https://github.com/Geoffrey-McGregor/CHMC-Codes}{https://github.com/Geoffrey-McGregor/CHMC-Codes}


\acknow{The authors acknowledge support from Natural Sciences and Engineering Research Council of Canada Discovery Grant (RGPIN-2019-07286) and the University of Northern British Columbia, where this work was initiated. A.T.S.W. acknowledges support from the University of California, Merced.}

\showacknow 

\bibliography{refs}

\includepdf[linktodoc=true,pages=-]{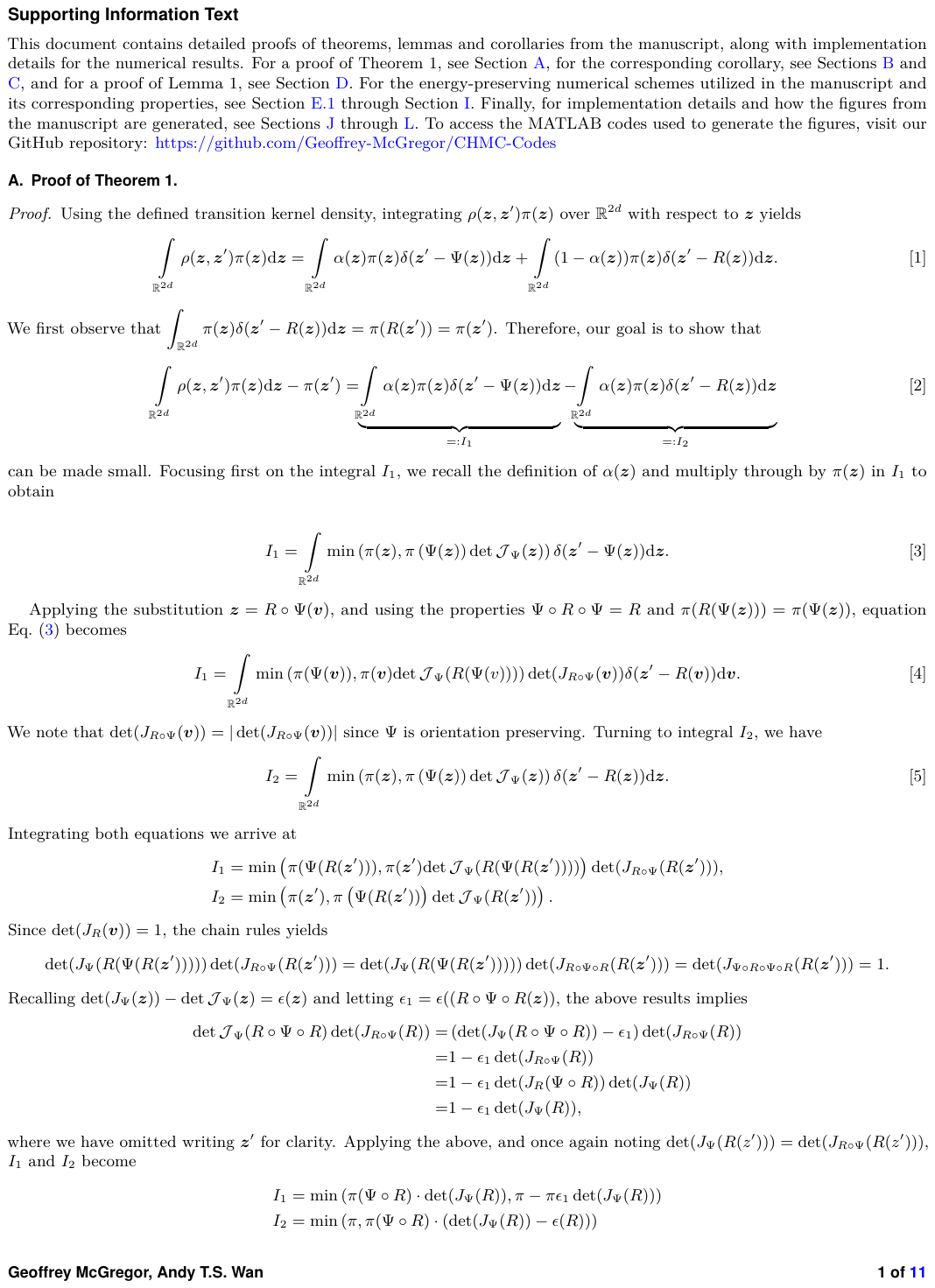}

\end{document}